\documentclass[12pt]{article}
\usepackage{latexsym, amsfonts, amsmath, amssymb, amscd, epsfig,extarrows}
by-\headheight \advance \topmargin by-\headsep 
\usepackage{graphicx}
\usepackage{amsmath,amssymb,amsthm,amsfonts,mathrsfs}
\usepackage{amssymb}
\newtheorem{thm}{Theorem}[section]

\newtheorem{lem}[thm]{Lemma}

\theoremstyle{definition}
\newtheorem{defn}{Definition}[section]
\theoremstyle{remark}
\newtheorem{rem}{Remark}[section]

\numberwithin{equation}{section}

\DeclareMathSymbol{\C}{\mathalpha}{AMSb}{"43}

\newcommand{\ef}{\eqref}

\newcommand{\fr}{\frac}

\newenvironment{pf}{{\noindent \it \bf Proof.\ }}{{\hfill$\Box$}\\}

\newcommand{\bsub}{\begin{subequations}}
	\newcommand{\esub}{\end{subequations}$\!$}

\begin{document}
\vspace{1cm}
\title{On the vanishing dissipation limit for the incompressible MHD equations on bounded domains
}

\author{Qin Duan$^{1,2}$
  Yuelong Xiao$^{1}$\thanks{\footnotesize {Corresponding author: xyl@xtu.edu.cn.} }
  and Zhouping Xin$^{3}$\\[1.5mm]
\footnotesize  { $^1$ Hunan Key Laboratory for Computation and Simulation in Science and Engineering, }\\[-2.8mm]
\footnotesize  {School of Mathematics and Computational Science, Xiangtan University}\\[-2.8mm]
\footnotesize  { $^2$ College of Mathematics and Statistics, Shenzhen University}\\[-2.8mm]
\footnotesize  { $^3$ The Institute of Mathematical Sciences, The Chinese University of Hong Kong}
}


\date{}




\maketitle

\begin{abstract}
In this paper, we investigate the solvability, regularity and the vanishing dissipation limit of solutions to  the three-dimensional viscous magneto-hydrodynamic (MHD) equations in  bounded domains. On the boundary, the velocity field fulfills a Navier-slip condition, while the magnetic field satisfies  the insulating condition. It is  shown that  the initial-boundary problem has a global weak solution for a general smooth domain. More importantly, for a flat domain, we establish the uniform local well-posedness of the strong solution with higher order uniform regularity and the asymptotic convergence with a rate  to the solution of the ideal MHD as the dissipation tends to zero.

\end{abstract}

\section{Introduction}
  \quad \quad Let $\Omega\subset\mathbb{R}^3$ be a bounded smooth domain with boundary $\partial\Omega$  and $n$ be the outward normal vector on $\partial\Omega$. The 3-dimentional (3-D) viscous magnetic-hydrodynamic equations (MHD) can be written as
\begin{equation}\label{1.1}\left\lbrace
\begin{split}
&\partial_tu-\nu\Delta u+(\nabla\times u)\times u+B\times(\nabla\times B)+\nabla P=0\quad \text{in}\  \Omega,\\
&\nabla\cdot u=0\quad\text{in}\  \Omega,\\
&\partial_tB-\mu\Delta B=\nabla\times(u\times B)\quad \text{in}\  \Omega,\\
&\nabla\cdot B=0  \quad \text{in}\  \Omega,
\end{split}\right.
\end{equation}
while $u$ and $B$ are the velocity field and magnetic field respectively, $P$ is the pressure, $\nu$ is the viscosity coefficient and $\mu$ is the magnetic diffusion coefficient, $\nabla\cdot$ and $\nabla\times$ denote the $\text{div}$ and $\text{curl}$ operators respectively. \ef{1.1} is supplemented with the following initial data
\begin{equation}\label{id}
(u,B)(x,t=0)=(u_0,B_0)(x),\ \quad x\in \Omega.
\end{equation}
On the boundary, the velocity field is assumed to satisfy the following Navier-slip condition (\cite{XX2,XXW}):
\begin{equation}\label{1.2}
u\cdot n=0,\quad \nabla\times u\times n=0\quad \text{on}\ \partial\Omega,
\end{equation}
and the magnetic field satisfies the insulating boundary condition (\cite{Gun,JeanF06})
\begin{equation}\label{1.3}
B\times n=0\quad \text{on}\ \partial\Omega,
\end{equation}
where $n$ is the unit outnormal of $\partial\Omega$.

The aim of this paper is to study the solvability, regularity and the asymptotic behavior as the dissipations vanish ($\nu\rightarrow 0^+$ and $\mu\rightarrow 0^+$) of the solution to the initial-boundary-value problem (IBVP) $\ef{1.1}-\ef{1.3}$. In particular, we are concerned with the uniform (with respect to the dissipation $\nu$ and $\mu$) well-posedness of the strong solution to the IBVP \ef{1.1}-\ef{1.3} and whether there are strong boundary layers near the physical boundary.

This study is motivated strongly by the vanishing viscosity limit problem for the 
incompressible Navier-Stokes equations, which is a classical problem in the mathematical theory of fluid dynamics and has been studied extensively in the case without physical boundaries (see \cite{Co1,Co,Ka,Ka1,Mas} for instance). Yet the vanishing viscosity problem for the Navier-Stokes equation becomes more complicated and challenging in the presence of physical boundaries despite its fundamental importance both physically and theoretically in understanding the boundary layer behavior of viscous flows for large Renolds number. Substantial difficulties arise due to the appearance of strong boundary layers in the case that the commonly used no-slip condition is imposed on the boundary, which makes it extremely difficult to justify the well-known Prandtl's boundary layer theory with few notable exception (\cite{Ka1,Mae,Sm1,Sm2,Wa}). On the other hand, if a Navier-slip condition, such as \ef{1.2}, is imposed on the boundary, then vorticities created near the boundary can be controlled so the possible viscous boundary layers are weak, which makes it possible to obtain strong convergence of the viscous fluid to the ideal one, as shown in many recent works, see (\cite{BF}-\cite{CMR}, \cite{IfP,IfS,Ke1,Ke2,Mas1,WWX,WXZ,XX,XX1,XX3,XX4}) and the references therein. In particular, for a flat domain $\Omega\subset\mathbb{R}^3$ with the slip boundary condition $\ef{1.2}$, Xiao and Xin (\cite{XX}) introduced an argument to obtain the uniform $H^3(\Omega)$ estimates on the solution to the 3-D Navier-Stokes equations and further the convergence with an optimal rate to the solution to the incompressible Euler system as the viscosity tends to zero. One is also referred to (\cite{BF,BF1,WXZ,XX,XX1}) for further studies in these directions.

As for the hydrodynamic case, the vanishing dissipation problem for the MHD equations in the presence of physical boundaries is important both in theory and applications, and is also a challenging topic. Due to the great complexities of the MHD system, there have been less results in the rigorous mathematical treatment of boundary layers for the MHD compared to the hydrodynamics until recently. Even for fixed positive dissipation, i.e., $\nu>0$, $\mu>0$, the well-posedness of initial boundary value problem (IBVP) for MHD has been developed mostly for the case that the velocity field $u$ satisfies the no-slip condition and while the magnetic field satisfies the perfect conducting condition, i.e.,
\begin{equation}\label{noslip}
u=0, \ \ \text{on}\ \ \partial\Omega,
\end{equation}
\begin{equation}\label{pc}
B\cdot n=0, \ \ \nabla\times B\times n=0 \ \ \ \text{on}\ \ \partial\Omega,
\end{equation}
see \cite{DL,JeanF06,ST} and the references therein. Indeed, it is  a subtle problem to prescribe boundary conditions for the magnetic field $B$ mathematically since the magnetic field $B$ satisfies a system of parabolic equation up to leading order with the additional divergence free constraint, i.e.,
\begin{equation}\label{B}\left\lbrace
\begin{split}
&\partial_tB-\mu\Delta B=\nabla\times(u\times B)\quad \text{in}\  \Omega,\\
&\nabla\cdot B=0 \quad \text{in}\  \Omega,
\end{split}\right.
\end{equation}
so the standard parabolic boundary condition for $B$ on $\partial\Omega$, such as the homogeneous Dirichlet boundary condition, may lead to an overdetermined problem. In fact, as far as we are aware, all the known well-posedness theory for the initial boundary value problem for the unsteady MHD, except \cite{HHM}, deals with boundary conditions as either \ef{noslip}, \ef{pc} (\cite{DL,JeanF06,ST}), or \ef{1.2} and \ef{pc}, (\cite{XXW}), or
\begin{equation}\label{B}
B\cdot n=0,\quad (\nabla\times B)\cdot n=0, \quad \Delta B\cdot n=0,\quad \text{on}\ \partial \Omega,
\end{equation}
see \cite{StoG}. In all these three cases, the boundary conditions for the magnetic field are characterized as that  with these boundary conditions for $B$, the corresponding Stokes and Laplacian operators are identical, which is not true in general for the Dirichlet boundary conditions, see
(\cite{Co,PLL96}). Thus we do not understand the argument in \cite{HHM}. For the stationary MHD, the well-posedness has been well established in the case that the velocity field satisfies the no-slip (Dirichlet) condition and the magnetic field satisfies either the perfectly conducting condition \ef{pc}, or the perfect insulating condition \ef{1.3}, see
\cite{Gun} and \cite{JeanF06} and references therein for generalizations. For more discussions on boundary conditions for MHD either mathematically or physically, we refer to (\cite{HY,Ja,JeanF06,she,StoG}) and references therein. To our knowledge, the well-posedness of the IBVP, \ef{1.1}-\ef{1.3} has not been considered before. This will be one of the consequences of the studies in this paper.

In the case that the dissipation coefficients are positive but can be arbitrarily small, i.e. $0<\nu, \mu\ll 1$, the uniform well-posedness and asymptotic behavior of solutions to the initial-boundary value problem for the MHD are difficult to study due to the possible appearance of boundary layers. The first rigorous result along this line is due to Xiao-Xin-Wu \cite{XXW} in  where they treat 3-dimensional flat domains with the slip boundary conditions, $\ef{1.2}$, for the velocity field and the perfect conducting boundary conditions, \ef{pc}, for the magnetic field, get the local uniform well-posedness of the solutions to the IBVP for the MHD with uniform $H^3(\Omega)-$ estimate independent of the dissipations, and finally obtain the asymptotic convergence with an optimal rate as the dissipations tend to zero, just as the corresponding results for the hydrodynamic case established by Xiao-Xin in \cite{XX}. However, it should be noted that in the theory of \cite{XXW}, there are no strong boundary layers due to the slip boundary condition as \ef{1.2}. In the case that the velocity field satisfies the no-slip boundary condition \ef{noslip}, strong boundary layers are expected and it is much harder to study the zero dissipation
limit problem. However, recently, for the case that the velocity field satisfies the no-slip conditions \ef{noslip} and the magnetic field satisfies \ef{pc} for the domain which is a half plane, Liu-Xie-Yang \cite{LXY} solved successfully the zero dissipation limit problem for the MHD system (locally) provided the magnetic field is non-trivial by solving the Prandtl's boundary layer equations for MHD first. Note that the theory in \cite{LXY} depends crucially on the presence of the magnetic field so that the results do not hold for the hydro-dynamic equations.

The main results in the paper concern with solutions to the MHD equation \ef{1.1} on a smooth bounded domain, subject to the slip boundary condition \ef{1.2} for the velocity field $u$ and the perfect insulating boundary condition  \ef{1.3} for the magnetic field $B$ with suitable initial condition \ef{id}. First, for fixed positive viscosity and magnetic diffusion coefficient and general smooth bounded domains, we show that there exists a global in time weak solution to the IBVP \ef{1.1}-\ef{1.3} for general initial data, and such a weak solution becomes the unique strong solution for short time for general regular initial data. As for the Navier-Stokes equations, such a local strong solution can be extended globally in time for suitable small initial data (depending on the dissipation coefficients). Second, for flat domains, $\Omega =\mathbb{T}^2\times(0,1)$, with $\mathbb{T}^2$ being the torus, we will derive the uniform (independent of the dissipation coefficients) $H^3(\Omega)$ estimates on the solutions (and thus the uniform well-posedness), which then enable us to prove the convergence of the solutions of IBVP \ef{1.1}-\ef{1.3} to  the solution of  the ideal MHD equations as $\nu$ and $\mu$ tend to zero. Furthermore, the optimal rate of convergence is obtained also, see Theorem \ref{uniform}- Theorem\ref{h4}. These are major results in this paper.

We now make some comments on the analysis of this paper. The existence of weak and strong solutions is proved by using the Galerkin method based on eigenvalue problems for the corresponding Stokes operators associated with boundary conditions \ef{1.2} and \ef{1.3} respectively, see Lemma 2.3 and Lemma 2.4. For the uniform well-posedness and vanishing dissipation limit, the basic idea is that the vorticity created near the boundary due to the slip condition $\ef{1.2}$ is weak and can be controlled as in the case for the Navier-Stokes equation investigated by Xiao-Xin in \cite{XX} provided that the effects of the magnetic field can be taken care of. Thus the results and the approach are similar to the case in \cite{XXW} where the velocity field and magnetic field satisfy the boundary condition, \ef{1.2} and \ef{pc} respectively. However, one cannot apply the analysis in \cite{XXW} into the case here directly. This is due to that the boundary conditions \ef{1.2} and \ef{pc} have the same structure so that one can use the argument of \cite{XX} to the system satisfied by velocity field and the magnetic field separately to derive the uniform $H^3(\Omega)$ estimates in \cite{XXW}. However, in our case, the perfect insulating conditions \ef{1.3} for the magnetic field are completely different from the slip condition \ef{1.2} for the velocity field. Thus to obtain the uniform high order estimates on the solutions to the IBVP \ef{1.1}-\ef{1.3} by energy method, one needs to use the equations for the velocity field  and the magnetic field and the corresponding boundary condition simultaneously. Indeed, it turns out that to derive the a priori estimates on the $H^3(\Omega)$ norm of the solutions to \ef{1.1}-\ef{1.3} by suitable energy method, it is crucial to have
\begin{equation}\label{ub1}
(\nabla\times)^3u\times n=0,\ \ (\nabla\times)^3(B\times u)\times n=0\ \ \text{on}\ \ \partial\Omega,
\end{equation}
which do not come from the boundary conditions \ef{1.2}, \ef{1.3} directly, but hold for strong solutions to the IBVP \ef{1.1}-\ef{1.3}. Indeed, $(\nabla\times)^3u\times n=0$ follows from the equations for the velocity field in \ef{1.1} and the boundary conditions \ef{1.2} and \ef{1.3} (see Lemma 5.1) provided that
\begin{equation}\label{ub2}
\nabla\times(B\times(\nabla\times B))\times n|_{\partial \Omega}=0.
\end{equation}
To show \ef{ub2}, one needs to use that the equations for the magnetic field and Lemma 5.1 to obtain that $\triangle B\times n|_{\partial\Omega}=0$. Then \ef{ub2} follows from this and the boundary condition \ef{1.2} and \ef{1.3}, see Lemma 5.2. Thus $(\nabla\times)^3u\times n|_{\partial \Omega}=0$ follows, which, together with \ef{1.2} and \ef{1.3}, implies $(\nabla\times)^3(B\times u)\times n|_{\partial\Omega}=0$, Lemma 5.3. As  consequences of \ef{ub1} and the observation of some cancellations of nonlinear terms (see \ef{*}), one can use the energy method to derive the desired uniform $H^3(\Omega)$ estimates and the convergence theory.

The rest of the paper is organized as follows. In the next section, we introduce the functional spaces associated with the boundary conditions \ef{1.2} and \ef{1.3} and study the eigenvalue problem for the corresponding Stokes operators associated with the boundary condition \ef{1.2} and \ef{1.3} respectively. The global existence of  a weak solution to the IBVP, \ef{1.1}-\ef{1.3}, is obtained in section 3, while the well-posedness of local (in time) strong solutions and regularity properties are established in section 4. Finally, we establish the uniform well-posedness and vanishing dissipation limit results in section 5, which are the main results of this paper.

\section{Preliminaries}
  \qquad In this section, we will introduce some basic tools that will be used for  later analysis.
Let $\Omega\subset R^3$ be a smooth domain and denote $H^s(\Omega)$ with $s\geq 0$ the standard Sobolev spaces and $H^{-s}(\Omega)$ with $s\geq 0$ denotes the dual of $H^s_0(\Omega)$. It is well known that:
\begin{lem}\label{us1}
	 Let $s\geq 0$ be an integer. Let $u\in H^s$ be a vector-valued function. Then
\begin{equation}\label{2.1}
\|u\|_s\leq  C(\|\nabla\times u\|_{s-1}+\|\nabla\cdot u\|_{s-1}+|u\cdot n|_{s-\fr{1}{2}}+\|u\|_{s-1}).
\end{equation}
\end{lem}

\begin{lem}\label{us2}Let $s\geq 0$ be an integer. Let $u\in H^s(\Omega)$. Then
\begin{equation}\label{2.2}
\|u\|_s\leq  C(\|\nabla\times u\|_{s-1}+\|\nabla\cdot u\|_{s-1}+|u\times n|_{s-\fr{1}{2}}+\|u\|_{s-1}).
\end{equation}
\end{lem}
See \cite{XX,XXW} and the references therein.

Let
\begin{equation*}
X=\{u\in L^2;\nabla\cdot u=0,u\cdot n=0\ \ \text{on}\ \partial\Omega\},
\end{equation*}
be the Hilbert space with the $L^2$ inner product, and let
\begin{equation*}
V=H^1\cap X\subset X,
\end{equation*}
\begin{equation*}
W=\{u\in X\cap H^2;(\nabla\times u)\times n=0 \ \text{on}\ \partial\Omega\}\subset X.
\end{equation*}
\begin{equation*}
\hat{X}=\{x\in L^2;\nabla\cdot B=0\},
\end{equation*}
\begin{equation*}
\hat{V}=\{x\in H^1;\nabla\cdot B=0,B\times n=0\ \ \text{on}\ \partial\Omega\},
\end{equation*}
\begin{equation*}
\hat{W}=H^2\cap \hat{V}\subset \hat{V}.
\end{equation*}
It follows from \ef{2.1} and \ef{2.2} that for any $u\in H^s(\Omega)\cap (V\cup \hat{V})$ , it holds that
\begin{equation}\label{2.3}
\|u\|_s\cong \|u\| + \|(\nabla\times)^s u\|.
\end{equation}

\begin{lem}\label{st} The Stokes operator $A=I-P\triangle = I - \Delta$ with its domain by $D(A)=W\subset V$ satisfying
\begin{equation*}
(Au,v)=a(u,v)\equiv(u,v)+\int_\Omega\nabla\times u\cdot\nabla\times vdx
\end{equation*}
is a self-adjoint and positive defined operator, with its inverse being compact. Consequently, its countable eigenvalues can be listed as
\begin{equation*}
0\leq\lambda_1\leq\lambda_2\leq \cdots\longrightarrow\infty
\end{equation*}
and the corresponding eigenvectors $\{e_j\}\subset W\cap C^\infty$,i.e.,
\[   Ae_j = (1+\lambda_j) e_j,\ or\  -\Delta e_j = \lambda_j e_j, \]
which form an orthogonal complete basis of $X$.
\end{lem}
For the details, refer to \cite{XX2}.

\begin{rem}
This Stokes operator $P\Delta = \Delta$ with the domain $D(A)=W$ follows from $\Delta u\cdot n|_{\partial\Omega}= 0$ for $u\in W$.
See \cite{XX2} for the details.
\end{rem}
Set
\begin{equation*}
\hat{e}_j=\frac{\nabla\times e_j}{\|\nabla\times e_j\|},\ \ j\in N.
\end{equation*}
It then follows that 
\begin{lem} $\{\hat{e}_j\}, j\in N$ form an orthogonal complete basis for $\hat{X}$. The bilinear form
\[ \hat{a}(u,v) = (u,v) + \int_\Omega (\nabla\times u)\cdot(\nabla\times v),\ \ \ \ \ u,v \in \hat{V}\]
with the domain $D(\hat{a})=\hat{V}$ is positive and closed, $\hat{a}$ is densely defined in $\hat{X}$. $\hat{A}= I - \Delta$ is the self-adjoint extension of the bilinear form $\hat{a}(u,v)$ with domain
$D(\hat{A})=\hat{W}$. Also $\hat{e}_j$ is the eigenvector of $\hat{A}$ respect to the eigenvalue $\lambda_j$.
\end{lem}
\pf It is clear that $\{\hat{e}_j\}, j\in N$ is an orthogonal series in $\hat{X}$.
Let $x\in \hat{X}$ and $(x,\nabla\times e_j)=0, \forall j\in N$.
Since $\nabla\cdot x = 0$ and $x\in L^2(\Omega)$,
one has
\[   x = \nabla\times u  \]
for some $u\in V$. It follows from
\[   0=(x,\nabla\times e_j) = (\nabla\times u,\nabla\times e_j) = (u,-\Delta e_j) = \lambda_j(u,e_j), \quad\forall j\in N, \]
that $u = 0$, and then $x = 0$. Since $\hat{X}$ is a Hilbert space, it follows that $\{\hat{e}_j\}, j\in N$ form an orthogonal complete basis of $\hat{X}$.
Clearly, $\hat{a}(u,v)$ is positive. Since $\hat{e}_j \subset D(\hat{a})$, for $j\in N$, then $\hat{a}$ is densely defined in $\hat{X}$.
From \ef{2.3}, one can define the bilinear form $\hat{a}(u,v)$ as a inner product on $\hat{V}$. By using of the trace theorem, and noting the continuity of the divergence operator,
it follows that $\hat{V}$ is closed in $H^1(\Omega)$, and then is a Hilbert space respect to the inner product. Hence, $\hat{a}(u,v)$ is closed.
It is clear that $\hat{e}_j$ also satisfies
\[    -\Delta \hat{e}_j = \lambda_j\hat{e}_j    \]
and the boundary condition
\[  \hat{e}_j\times n=0 \]
on the boundary.
We can define the operator $\hat{A}: \overline{span\{\hat{e}_j\}}^2\mapsto\hat{X}$ by
\begin{equation}
  \hat{A}(\sum u_j \hat{e}_j)= \sum u_j(1-\Delta)\hat{e}_j =  \sum u_j (1+\lambda_j)\hat{e}_j
  \end{equation}
for $u= \sum u_j \hat{e}_j\in \overline{span\{\hat{e}_j\}}^2$, where $\overline{span\{\hat{e}_j\}}^2$ is the closure of $\{\hat{e}_j\}$ in $H^2(\Omega)$.
It is easy to check that $\overline{span\{\hat{e}_j\}}^2 = \hat{W}$. Indeed, let $v \in \hat{W}$ and
\[   (\hat{e}_j,v) + (\Delta \hat{e}_j,\Delta v) = 0, \ \ \forall j.  \]
It follows that
\[   (1+\lambda_j^2)(\hat{e}_j,v) = 0,\ \ \ \forall j.  \]
Since $\hat{W}$ is a Hilbert space respect to the inner product
\[  (u,v)_{\hat{A}} = (u,v) + (\Delta u,\Delta v), \]
 so we have $v = 0$ and $\overline{span\{\hat{e}_j\}}^2 = \hat{W}$.

Now, let $f\in \hat{X}$, we can write
\[    f = \sum f_j \hat{e}_j.   \]
It follows that $u = \sum u_j\hat{e}_j\in \hat{W}$, with $(1+\lambda_j)u_j = f_j$ is a solution of $\hat{A}u= f$ in the sense that
\[     \hat{a}(u,v) = (f,v), \ \ \ \forall v\in \hat{V}.    \]
If $f= 0$, it follows that
\[   (u,\hat{e}_j) + \lambda_j(u,\hat{e}_j) = 0,  \]
by taking $v = \hat{e}_j$.
Then, $u = 0$, and we conclude that $\hat{A}:\hat{W}\mapsto \hat{X}$ is a isomorphism. The lemma is proved.

For $u\in W$, $B\in\hat{W}$, we  define
\begin{equation*}
H_1:W\times\hat{W}\longrightarrow X
\end{equation*}
by
\begin{equation*}
H_1(u,B)=(\nabla\times u)\times u+B\times(\nabla\times B)+\nabla p,
\end{equation*}
where $p$ satisfies
\begin{equation*}
\begin{split}
&\triangle p=-\nabla\cdot((\nabla\times u)\times u+B\times(\nabla\times B)),\\
&\nabla p\cdot n=-(\nabla\times u)\times u+B\times(\nabla\times B)\cdot n  \ \ \text{on} \ \ \partial \Omega
\end{split}
\end{equation*}
and
\begin{equation*}
H_2:W\times\hat{W}\longrightarrow \hat{X}
\end{equation*}
by
\begin{equation*}
H_2(u,B)=\nabla\times (B\times u).
\end{equation*}

\section{The weak solutions}
\qquad In this section, we will establish the global existence of weak solutions to the  systems $\ef{1.1}-\ef{1.3}$ by using of the method of Galerkin approximation. Here we consider a general smooth bounded domain in $R^3$ unless stated otherwise. We denote $H^*$ be the dual space of $H$.

\begin{defn}\label{weak} $(u,B)$ is called a weak solution of $\ef{1.1}-\ef{1.3}$ with the initial data $u_0\in X$, $B_0\in \hat{X}$ on the time interval $[0,T)$ if $u\in L^2(0,T;V)\cap C_w([0,T);X)$ and $B\in L^2(0,T;\hat{V})\cap C_w([0,T);\hat{X})$ satisfy $u'\in L^1(0,T;V^*)$, $B'\in L^1(0,T;\hat{V}^*)$ and
\begin{equation}\label{ws}
\begin{split}
&(u',\phi)+\nu(\nabla\times u,\nabla\times\phi)+(\nabla\times u\times u,\phi)+(B\times\nabla\times B,\phi)=0,\\
&(B',\psi)+\mu(\nabla\times B,\nabla\times\psi)+(u\cdot\nabla B-B\cdot\nabla u,\psi)=0
\end{split}
\end{equation}
for all $\phi\in V$, $\psi\in \hat{V}$ and for a.e. $t\in[0,T)$, and
\begin{equation*}
u(0)=u_0,\quad B(0)=B_0.
\end{equation*}
\end{defn}

\begin{thm}\label{weak-thm} Let $u_0\in X$, $B_0\in \hat{X}$. Let $T>0$. Then there exists at least one weak solution $(u,B)$ of $\ef{1.1}-\ef{1.3}$ on $[0,T)$ which satisfies the energy inequality
\begin{equation}\label{3.1}
\fr{d}{dt}(\|u\|^2+\|B\|^2)+2(\nu\|\nabla\times u\|^2+\mu\|\nabla\times B\|^2)\leq 0
\end{equation}
in the sense of the distribution.
\end{thm}
\pf This will be proved by using a Galerkin approximation based on eigenvalue problems in Lemma 2.2 and Lemma 2.3. Define
\begin{equation*}
u^m(t)=\sum^m_{j=1}u_j(t)e_j,\quad  B^m(t)=\sum^m_{j=1}b_j(t)\hat{e}_j,
\end{equation*}
where $u_j$ and $b_j$ for $j=1,2,\cdots,m$ solve the following ordinary differential equations
\begin{equation}\label{3.2}
u_j'(t)+\nu\lambda_j u_j(t)+g^1_j(U)=0,
\end{equation}
\begin{equation}\label{3.3}
b_j'(t)+\mu\lambda_j b_j(t)+g^2_j(U)=0,
\end{equation}
\begin{equation}\label{3.4}
u_j(0)=(u_0,e_j),\quad b_j(0)=(B_0,\hat{e_j}),
\end{equation}
with $U=(u_1,u_2,\cdots,u_m,B_1,B_2,\cdots B_m)$ and
\begin{equation*}
g^1_j(U)=(H_1(u^m,B^m),e_j),
\end{equation*}
\begin{equation*}
g^2_j(U)=(H_2(u^m,B^m), \hat{e}_j).
\end{equation*}
It follows from the Lipschitz continuity of $(g^k_j(U))$ that the initial value problem $\ef{3.2}-\ef{3.4}$ is locally well-posed on interval $[0,T)$ for some positive $T$. Consequently, for any $t\in [0,T)$, $(u^m,B^m)$ solves the following systems of equations
\begin{equation}\label{3.5}
(u^m)'-\nu\triangle u^m+P_mH_1(u^m,B^m)=0,
\end{equation}
\begin{equation}\label{3.6}
(B^m)'-\mu\triangle B^m+\hat{P}_mH_2(u^m,B^m)=0,
\end{equation}
\begin{equation}\label{3.7}
u^m(0)=P_mu_0,\quad B^m(0)=\hat{P}_mB_0
\end{equation}
where $P_m$ denotes the projection of $X$ onto the space spanned by $\{e_j\}^m_1$ and  $\hat{P}_m$ denotes the projection of $\hat{X}$ onto the space spanned by $\{\hat{e}_j\}^m_1$.

Multiplying  $u^m$ and $B^m$ on both side of $\ef{3.5}$ and $\ef{3.6}$ respectively and integrating by parts, lead to
\begin{equation*}
(P_mH_1(u^m,B^m),u^m)=\int_\Omega(B^m\times \nabla\times B^m)\cdot u^mdx,
\end{equation*}
\begin{equation*}\begin{split}
(\hat{P}_mH_2(u^m,B^m),B^m)&=\int_\Omega\nabla\times(B^m\times u^m)\cdot B^mdx.\\
&=-\int_\Omega (B^m\times\nabla\times B^m)\cdot u^mdx.
\end{split}
\end{equation*}
Adding them up, one obtains by a simple computation that
\begin{equation}\label{3.8}
\fr{d}{dt}(\|u^m\|^2+\|B^m\|^2)+2(\nu\|\nabla \times u^m\|^2+\mu\|\nabla \times B^m\|^2)=0,
\end{equation}
which implies that
\begin{equation}\label{3.9}
u^m \ \text{is bounded in }\ L^\infty(0,T;X),\   B^m\ \text{is bounded in }\ L^\infty(0,T;\hat{X}),\ \text{uniformly for m}.
\end{equation}
\begin{equation}\label{3.10}
 u^m \ \text{is bounded in }\ L^2(0,T;V), \ \  B^m \ \text{is bounded in }\ L^2(0,T;\hat{V}), \ \text{uniformly for m}.
\end{equation}
Note that for $\phi\in V$ and $\psi\in \hat{V}$, it holds that
\begin{equation*}
|(-\triangle u^m,\phi)|=|(\nabla\times u^m,\nabla\times \phi)|,
\end{equation*}
\begin{equation*}
|(-\triangle B^m,\psi)|=|(\nabla\times B^m,\nabla\times \psi)|.
\end{equation*}
Therefore,
\begin{equation*}
\{-\triangle u^m\} \ \text{is bounded in }\ L^2(0,T;V^*),
\end{equation*}
and
\begin{equation*}
\{-\triangle B^m\} \ \text{is bounded in }\ L^2(0,T;\hat{V}^*).
\end{equation*}
On the other hand, it follows from Sobolev inequalities that for any  $\phi\in V$,
\begin{equation*}\begin{split}
&|(P_mH_1(u^m,B^m),\phi)|=|(H_1(u^m,B^m),\phi_m)| \\
&\leq C(\|u^m\|^{\fr{1}{2}}\|u^m\|^{\fr{3}{2}}_1+\|B^m\|^{\fr{1}{2}}\|B^m\|^{\fr{3}{2}}_1)\|\phi_m\|_1,
\end{split}
\end{equation*}
where $\phi_m=P_m\phi$. From the uniform bound of $\ef{3.8}$, we have
\begin{equation*}
\{H_1(u^m,B^m)\}\ \text{is bounded in }\ L^\fr{4}{3}(0,T;V^*).
\end{equation*}
Similarly,
\begin{equation*}
\{H_2(u^m,B^m)\}\ \text{is bounded in }\  L^\fr{4}{3}(0,T;\hat{V}^*).
\end{equation*}
Hence,
\begin{equation*}
(u^m)' \ \text{is bounded in }\ L^\fr{4}{3}(0,T;V^*),
\end{equation*}
\begin{equation*}
(B^m)' \ \text{is bounded in }\ L^\fr{4}{3}(0,T;\hat{V}^*).
\end{equation*}
Now we can use the similar arguments in Constantin and Foias in \cite{Co} to complete the proof of Theorem \ref{weak-thm}. The details are  omitted.

\section { The strong solutions}
\qquad In this section, we will study the local well-posedness of strong solutions of $\ef{1.1}-\ef{1.3}$ with the corresponding initial data $u_0\in V$ and $B_0\in \hat{V}$.
\begin{defn}\label{strong}
 $(u,B)$ is called a strong solution of $\ef{1.1}-\ef{1.3}$ with the initial data $u_0\in V$, $B_0\in \hat{V}$ on the time interval $[0,T)$ if $u\in L^2(0,T;W)\cap C([0,T);V)$ and $B\in L^2(0,T;\hat{W})\cap C([0,T);\hat{V})$ satisfy $u'\in L^2(0,T;X)$, $B'\in L^2(0,T;\hat{X})$ and
\begin{equation}\label{ss2}
\begin{split}
&(\partial_tu,\varphi)-\nu(\Delta u,\varphi) +((\nabla\times u)\times u+B\times(\nabla\times B),\varphi) =0\\
&(\partial_t B,\psi)-\mu(\Delta B,\psi) = (\nabla\times(u\times B),\psi)
\end{split}
\end{equation}
for any $(\varphi,\psi)\in X\times\hat{X}$
 and for a.e. $t\in[0,T)$, and
\begin{equation*}
u(0)=u_0,\quad B(0)=B_0.
\end{equation*}	
\end{defn}
\begin{thm}\label{str} Let $u_0\in V$, $B_0\in \hat{V}$. Then there exists a time $T^*>0$ depending only on $\nu,\mu$ and the $H^1-$ norm of $(u_0,B_0)$ such that $\ef{1.1}-\ef{1.3}$ has a unique strong solution $(u,B)$ on $[0,T^*)$ satisfying
the energy equation
\begin{equation}\label{ee}
\begin{split}
&\fr{d}{dt}(\|\nabla\times u\|^2+\|\nabla\times B\|^2)+2(\nu\|\triangle u\|^2+\mu\|\triangle B\|^2)\\
&+2(\nabla\times H_1(u,B),\nabla\times u)+2(H_2(u,B),-\triangle B)=0.
\end{split}
\end{equation}

\end{thm}
\pf Taking the curl of $\ef{3.5}$ yields
\begin{equation}\label{4.1}
(\nabla\times u^m)'-\nu\triangle(\nabla\times u^m)+\sum_{j=1}^m g^1_j\nabla\times e_j=0.
\end{equation}
It follows from $\ef{3.6}$ that
\begin{equation}\label{4.2}
(B^m)'-\mu\triangle B^m+\sum^m_{j=1}g^2_j\hat{e}_j=0.
\end{equation}
Taking the inner product $(\ef{4.1},\nabla\times u^m)+(\ef{4.2},-\triangle B^m)$, noting that
\begin{equation*}
(\nabla\times e_i,\nabla\times e_j)=\lambda_j(e_i,e_j),
\end{equation*}
one can get
\begin{equation}\label{4.4}
\begin{split}
&\fr{d}{dt}(\|\nabla\times u^m\|^2+\|\nabla\times B^m\|^2)+2(\nu\|(\nabla\times)^2 u^m\|^2+\mu\|(\nabla\times)^2 B^m\|^2)\\
&+2((\nabla\times H_1(u^m,B^m),\nabla\times u^m)+(H_2(u^m,B^m),-\triangle B^m))=0.
\end{split}
\end{equation}
Since
\begin{equation*}
|(\nabla\times H_1(u^m,B^m),\nabla\times u^m)|\leq C(\|B^m\|_\infty\|B^m\|_1\|(\nabla\times)^2u^m\|+\|u^m\|_\infty\|u^m\|_1\|(\nabla\times)^2u^m\|),
\end{equation*}
and
\begin{equation*}
|( H_2(u^m,B^m),-\triangle B^m)|\leq C(\|u^m\|_\infty+\|B^m\|_\infty)(\|u^m\|_1+\|B^m\|_1)\|(\nabla\times)^2B^m\|.
\end{equation*}
By applying the Agmon inequality
\begin{equation*}
\|f\|_\infty^2\leq \|f\|_1\|f\|_2,\quad \forall f\in H^2,
\end{equation*}
the equivalent of the norms in (\ref{2.3}), and the standard interpolation inequalities, one can get
\begin{equation*}
\begin{split}
&|(\nabla\times H_1(u^m,B^m),\nabla\times u^m)|\leq C(\|u^m\|_1^{\fr{3}{2}}+\|B^m\|_1^{\fr{3}{2}})(\|\Delta u^m\|^{\fr{1}{2}}+\|\triangle B^m\|^{\fr{1}{2}})\|\triangle u^m\|,\\
&|(H_2(u^m,B^m),-\triangle B^m)|\leq C(\|u^m\|_1^{\fr{3}{2}}+\|B^m\|_1^{\fr{3}{2}})(1 + \|\Delta u^m\|^{\fr{1}{2}}+\|\triangle B^m\|^{\fr{1}{2}})\|\triangle B^m\|.
\end{split}
\end{equation*}
Combining with the $\ef{4.4}$, it holds
\begin{equation}\label{global}
\fr{d}{dt}(\|\nabla\times u^m\|^2+\|\nabla\times B^m\|^2)+\nu\|(\nabla\times)^2u^m\|^2+\mu\|(\nabla\times)^2B^m\|^2\leq C(1+ \|u^m\|_1+\|B^m\|_1)^6,
\end{equation}
where $C$ depends on $\nu$ and $\mu$.
Combining this with the energy inequality $\ef{3.1}$ shows that there is a time $T^*>0$ such that, for any fixed $T\in (0,T^*)$,
\begin{equation*}
\begin{split}
&(u^m,B^m)\  \text {is bounded in}\  L^\infty(0,T;H^1),\\
&(u^m,B^m)\  \text {is bounded in}\  L^2(0,T;H^2).
\end{split}
\end{equation*}
It then follows from
\begin{equation}\label{4.6}
\|P_m(u\times v)\|\leq \|u\times v\|\leq \|u\|_\infty\|v\|,
\end{equation}
\begin{equation}\label{4.6}
\|\hat{P}_mB\|\leq \|B\|,
\end{equation}
the definition of $H_1$ and $H_2$, and \ef{3.5}-\ef{3.7} that
\begin{equation}\label{4.7}
\{(u^m)'\},\{(B^m)'\} \ \text {is bounded in}\  L^2(0,T;L^2).
\end{equation}
Due to the Aubin-Lions Lemma, one can find a sequence of $(u^m,B^m)$ still denoted by $(u^m,B^m)$ and $(u,B)$ such that
\begin{equation*}
\begin{split}
&u^m\longrightarrow u\ \ \text{in}\  L^\infty(0,T;V) \ \ \text{weak-star},\quad B^m\longrightarrow B \ \text{in}\  L^\infty(0,T;\hat{V}) \ \text{weak-star},\\
&u^m\longrightarrow u \ \ \text{in}\  L^2(0,T;W) \ \text{weakly},\quad B^m\longrightarrow B \ \text{in}\  L^2(0,T;\hat{W}) \ \text{weakly},\\
&u^m\longrightarrow u\ \ \text{in}\  L^2(0,T;V) \ \ \text{strongly},\quad B^m\longrightarrow B \ \text{in}\  L^2(0,T;\hat{V}) \ \text{strongly}.\\
\end{split}
\end{equation*}
Passing to the limit, we find that $u\in L^\infty(0,T;V)\cap L^2(0,T;W)$, $B\in L^\infty(0,T;\hat{V})\cap L^2(0,T;\hat{W})$ satisfy
\begin{equation}\label{ss1}
\begin{split}
&(\partial_tu,\varphi)-\nu(\Delta u,\varphi) +((\nabla\times u)\times u+B\times(\nabla\times B),\varphi) =0\\
&(\partial_t B,\psi)-\mu(\Delta B,\psi) = (\nabla\times(u\times B),\psi)
\end{split}
\end{equation}
for any $(\varphi,\psi)\in X\times\hat{X}$.

From $\ef{4.7}$,
we find that $(u',B')\in L^2(0,T;X)\times L^2(0,T;\hat{X})$. Then we deduce that $(u,B)\in C([0,T];V)\times C([0,T];\hat{V})$, $u(0)=u_0$, $B(0)=B_0$ and the energy equation $\ef{ee}$ holds. 

Now we will show that the uniqueness of strong solutions by using of the standard procedure. Let $(u,B)$ and $(u^0,B^0)$ be two strong solutions to \ef{1.1}-\ef{1.3}. Set $\bar{u}=u-u^0$, and $\bar{B}=B-B^0$, then
\begin{equation}\label{4.8}
\partial_t\bar{u}-\nu\triangle \bar{u}+H_1(u,B)-H_1(u^0,B^0)=0,
\end{equation}
\begin{equation}\label{4.9}
\partial_t\bar{B}-\mu\triangle \bar{B}+H_2(u,B)-H_2(u^0,B^0)=0.
\end{equation}
Taking the inner products with $\bar{u}$ in $\ef{4.8}$, and  $\bar{B}$ in $\ef{4.9}$, integrating by parts, from the boundary condition $\bar{u}\cdot n=0$, $\bar{B}\times n=0$, one gets
\begin{equation*}
\begin{split}
&\fr{d}{dt}(\|\bar{u}\|^2+\|\bar{B}\|^2)+2\nu\|\nabla\times\bar{u}\|^2+2\mu\|\nabla\times\bar{B}\|^2\\
&+(H_1(u,B)-H_1(u^0,B^0),\bar{u})+(H_2(u,B)-H_2(u^0,B^0),\bar{B})=0.\\
\end{split}
\end{equation*}
On the other hand, simple manipulations show that 
\begin{equation}\label{h1h2}
\begin{split}
&(H_1(u,B)-H_1(u^0,B^0),\bar{u})+(H_2(u,B)-H_2(u^0,B^0),\bar{B})\\
&=(\nabla\times\bar{u}\times u,\bar{u})
+(B\times\nabla\times\bar{B},\bar{u})\\
&+(\bar{B}\times\nabla\times B^0,\bar{u})+\big(\nabla\times(B\times u)-\nabla\times(B^0\times u^0),\bar{B}\big)\\
&=(\nabla\times\bar{u}\times u,\bar{u})+(B\times\nabla\times\bar{B},\bar{u})+(\nabla\times(\bar{u}\times\bar{B}),B^0)\\
&+\big((B\times \bar{u},\nabla\times \bar{B})+(\bar{B}\times u^0,\nabla\times \bar{B})\big),\\
&=(\nabla\times\bar{u}\times u,\bar{u})+(\nabla\times(\bar{u}\times\bar{B}),B^0)+(\bar{B}\times u^0,\nabla\times \bar{B})),\\
\end{split}
\end{equation}
here one has used integration by parts, the boundary condition $\bar{B}\times n=0$, $B^0\times n=0$ and
\begin{equation*}
(\bar{B}\times\nabla\times B^0,\bar{u})=(\nabla\times(\bar{u}\times\bar{B}),B^0).
\end{equation*}
 Combining them together, yields
\begin{equation*}
\begin{split}
&\fr{d}{dt}(\|\bar{u}\|^2+\|\bar{B}\|^2)+2\nu\|\nabla\times\bar{u}\|^2+2\mu\|\nabla\times\bar{B}\|^2\\
&\leq \nu\|\nabla\times\bar{u}\|^2+\mu\|\nabla\times\bar{B}\|^2+C(\|u\|_2^2+\|B\|_2^2+\|u^0\|_2^2+\|B^0\|_2^2)(\|\bar{u}\|^2+\|\bar{B}\|^2)
\end{split}
\end{equation*}
on $[0,T^*)$. \\
Then $u=u^0$, $B=B^0$ follows from $\bar{u}(0)=\bar{B}(0)=0$ and the Gronwall's inequality.

\begin{rem}
	As usual, the local strong solution in Theorem \ref{str} can be extended globally in time if the initial data is suitably small. Indeed, note first that if $T^*$ be  the maximal time for the existence of strong solution in Theorem \ref{str} and $T^*<\infty$, then $\|(u,B)\|_1\longrightarrow 0$ as $t\longrightarrow T^*_{-}$. Next, it follows easily from  \ef{3.1} and \ef{ee}  that for $t\in(0, T^*)$,
\begin{equation}\label{ws}
\begin{split}
\fr{d}{dt}(\|u\|_1^2+\|B\|_1^2)&+\big(\mu-C(\|u\|_1+\|B\|_1)\big)\|u\|^2_2 \\
&+\big(\nu-C(\|u\|_1+\|B\|_1)\big)\| B\|^2_2
\leq 0.
\end{split}
\end{equation}
So the standard continuity argument shows that if 
\begin{equation*}
\|(u_0,B_0)\|_1\leq \delta(\mu,\nu),\ \ \text{then}\ \ T^*=+\infty,
\end{equation*}
so that the strong solution can be extended globally.

\end{rem}
\section {Uniform regularity and the inviscid limit}
 \qquad In last section, the existence time interval $[0,T^*)$ may depends on the parameters $\nu$ and $\mu$ for general domains. To study the uniform well-posedness independent of the viscosity and magnetic diffusion, we consider only the cubic domain $\Omega =\Omega_T=T^2\times(0,1)$, where $T^2$ is the torus. In this case, we can get the uniform well-posedness and convergence. We start with the following three elementary lemmas which will play an essential role in the proof of the Theorem \ref{uniform}.
\begin{lem}\label{bc1} Let $u,B$ be two smooth vectors satisfying
$\nabla\cdot u = 0$, $\nabla\cdot B = 0$
in $\Omega$ and
$u\cdot n = 0$, $(\nabla\times u)\times n = 0$; $B\times n = 0$ on the boundary,
that is
\begin{equation}\label{fu1}
\partial_3 u_1 = 0,\ \partial_3 u_2 = 0,\ u_3 = 0\ {\rm on} \ \partial\Omega,
\end{equation}
\begin{equation}\label{fb1}
	B_1 = 0,\ B_2 = 0,\ \partial_3 B_3 = 0\ {\rm on} \ \partial\Omega,
\end{equation}
where $\partial_3 B_3 = 0$ is from $\nabla\cdot B = 0$.
Then, it holds that
\begin{equation}\label{le1}
\begin{split}
&\nabla\times(B\times u)\times n = (u\cdot\nabla B - B\cdot\nabla u)\times n = 0\ {\rm on} \ \partial\Omega,\\
&\nabla\times(\nabla\times u\times u)\times n=0\ {\rm on} \ \partial\Omega.
\end{split}
\end{equation}
\end{lem}
\pf Indeed, direct calculations show that
\[ (u\cdot\nabla B )_1 = u_1\partial_1B_1 +  u_2\partial_2B_1 +  u_3\partial_3B_1= 0\ {\rm on} \ \partial\Omega,   \]
\[ (B\cdot\nabla u )_1 = B_1\partial_1u_1 +  B_2\partial_2u_1 +  B_3\partial_3u_1= 0\ {\rm on} \ \partial\Omega.   \]
Similarly
\[ (u\cdot\nabla B )_2 = u_1\partial_1B_2 +  u_2\partial_2B_2 + u_3\partial_3B_2= 0\ {\rm on} \ \partial\Omega,   \]
\[ (B\cdot\nabla u )_2 = B_1\partial_1u_2 +  B_2\partial_2u_2 +  B_3\partial_3u_2= 0\ {\rm on} \ \partial\Omega.   \]
Let $\omega=\nabla\times u$. Then
\begin{equation*}
\nabla\times(\nabla\times u\times u)\times n=\big((u\cdot \nabla \omega)-(\omega\cdot \nabla u)\big)\times n.
\end{equation*}
Since
\begin{equation*}
\omega_1=\omega_2=0 \ \ \text{on} \ \ \partial\Omega,
\end{equation*}
one gets easily  that
\begin{equation*}
\begin{split}
&(u\cdot\nabla\omega)_1=u_1\partial_1\omega_1+u_2\partial_2\omega_1+u_3\partial_3\omega_1=0 \ \ {\rm on} \ \partial\Omega,\\
&(u\cdot\nabla\omega)_2=u_1\partial_1\omega_2+u_2\partial_2\omega_2+u_3\partial_3\omega_2=0 \ \ {\rm on} \ \partial\Omega,\\
\end{split}
\end{equation*}
and
\begin{equation*}
\begin{split}
&(\omega\cdot\nabla u)_1=\omega_1\partial_1 u_1+\omega_2\partial_2 u_1+\omega_3\partial_3 u_1=0\ \ {\rm on} \ \partial\Omega,\\
&(\omega\cdot\nabla u)_2=\omega_1\partial_1 u_2+\omega_2\partial_2 u_2+\omega_3\partial_3 u_2=0\ \ {\rm on} \ \partial\Omega,
\end{split}
\end{equation*}
that is
\begin{equation*}
\nabla\times(\nabla\times u\times u)\times n=0\ \ {\rm on} \ \partial\Omega.
\end{equation*}

\begin{lem}\label{bc2} Let $u,B$ be the vectors in Lemma \ref{bc1}. Assume further that $\Delta B\times n = 0$ on the boundary,
that is
\begin{equation}\label{la}
	\partial^2_{33} B_1 = 0,\ \partial^2_{33} B_2 = 0\ {\rm on} \ \partial\Omega.
\end{equation}
Then, it holds that
\begin{equation}\label{le2}
	\nabla\times(\nabla\times B\times B)\times n = (B\cdot\nabla (\nabla\times B) - (\nabla\times B)\cdot\nabla B)\times n = 0\ {\rm on} \ \partial\Omega.
\end{equation}
\end{lem}
\pf It follows from \ef{la} and \ef{le2} that
\[   \nabla\times(\omega_B\times B)_1 = (B\cdot\nabla \omega_B - \omega_B\cdot\nabla B)_1,\]
\[  (B\cdot\nabla \omega_B)_1 = B_1\partial_1\omega_B^1+B_2\partial_2\omega_B^1 + B_3\partial_3\omega_B^1 \]
\[  = B_3\partial_3(\partial_2B_3-\partial_3B_2) =  B_3\partial^2_{32}B_3-B_3\partial^2_{33}B_2=0 \]
\[  (\omega_B\cdot\nabla B)_1 = \omega_B^1\partial_1B_1+\omega_B^2\partial_2B_1 + \omega_B^3\partial_3B_1 \]
\[  = (\partial_1B_2 - \partial_2B_1)\partial_3B_1 = 0,  \]
where $\omega_B = \nabla\times B$.

Similarly,
\[   \nabla\times(\omega_B\times B)_2 = (B\cdot\nabla \omega_B - \omega_B\cdot\nabla B)_2 = 0. \]
The lemma is proved.
\begin{lem}\label{bc3}  Let the assumptions in Lemma \ref{bc2} hold. Furthermore, $\omega=\nabla\times u$ satisfies
$\Delta \omega\times n = 0$ on the boundary. Then, it holds that
\begin{equation}\label{le3}
	(\nabla\times)^3(B\times u)\times n  = 0\ {\rm on} \ \partial\Omega.
\end{equation}
\end{lem}
\pf Direct computations yield
\begin{equation*}\begin{split}
&(\nabla\times)^3(B\times u)_1=-\Delta(u\cdot\nabla B-B\cdot\nabla u)_1\\
&=-(\partial^2_{11}+\partial^2_{22})(u\cdot\nabla B - B\cdot\nabla u)_1- \partial^2_{33}(u\cdot\nabla B - B\cdot\nabla u)_1 \\
&=-\partial^2_{33}(u\cdot\nabla B - B\cdot\nabla u)_1 \\
&  = -\partial^2_{33}(u_1\partial_1B_1+u_2\partial_2B_1 +u_3\partial_3B_1- B_1\partial_1u_1- B_2\partial_2u_1- B_3\partial_3u_1).
\end{split}
\end{equation*}

Note that
\begin{equation*}
\begin{split}
& \partial^2_{33}(u_1\partial_1B_1) =  \partial_3(\partial_3u_1\partial_1B_1 + u_1\partial^2_{31}B_1) \\
&  = \partial^2_{33}u_1\partial_1B_1+ 2\partial_3u_1\partial^2_{13}B_1 + u_1\partial_1\partial^2_{33}B_1=0,
\end{split}
\end{equation*}
\begin{equation*}
\begin{split}
& \partial^2_{33}(u_2\partial_2B_1) =\partial_3(\partial_3u_2\partial_2B_1 + u_2\partial^2_{32}B_1) \\
& = \partial^2_{33}u_2\partial_2B_1 + 2\partial_3u_2\partial^2_{23}B_1 + u_2\partial_2\partial^2_{33}B_1=0,
\end{split}
\end{equation*}
\begin{equation*}
\begin{split}
&\partial^2_{33}(u_3\partial_3B_1) =  \partial_3(\partial_3u_3\partial_3B_1 + u_3\partial^2_{33}B_1)\\
&= \partial^2_{33}u_3\partial_3B_1 + 2\partial_3u_3\partial^2_{33}B_1 + u_3\partial^3_{333}B_1=0
\end{split}
\end{equation*}
on the boundary, here one has used
\begin{equation*}
\begin{split}
&\partial^2_{33}u_3 = \Delta u_3 = -(\nabla\times \omega)_3\\
&= \partial_2\omega_1 - \partial_1\omega_2 =0.
\end{split}
\end{equation*}
Note also that
\begin{equation*}
\begin{split}
&\partial^2_{33}(B_1\partial_1u_1) =\partial^2_{33}B_1\partial_1u_1 + 2\partial_3B_1\partial_1\partial_{3}u_1 + B_1\partial_1\partial^2_{33}u_1=0,
\end{split}
\end{equation*}
\begin{equation*}
\begin{split}
&\partial^2_{33}(B_2\partial_2u_1)=\partial^2_{33}B_2\partial_2u_1 + 2\partial_3B_2\partial_{2}\partial_3u_1 + B_2\partial_2\partial^2_{33}u_1=0,
\end{split}
\end{equation*}
\begin{equation*}
\begin{split}
&\partial^2_{33}(B_3\partial_3u_1) = \partial^2_{33}B_3\partial_3u_1 + 2\partial_3B_3\partial^2_{33}u_1 + B_3\partial^3_{333}u_1=0
\end{split}
\end{equation*}
on the boundary, here one has used the fact
\begin{equation*}
\begin{split}
&\partial^3_{333}u_1 = \partial^2_{33}\partial_3u_1= (\partial^2_{33}\omega_2 + \partial_1\partial^2_{33}u_3) \\
&= \Delta\omega_2 + \partial_1(\Delta u_3-\partial^2_{11}u_3+\partial^2_{22}u_3) = 0
\end{split}
\end{equation*}
on the boundary, since
\begin{equation*}
-\Delta u_3 = -(\Delta u)_3 =(\nabla\times\omega)_3 = \partial_1\omega_2-\partial_2\omega_1
\end{equation*}
on the boundary.
Then, we conclude that
\begin{equation*}
(\nabla\times)^3(B\times u)_1 = 0
\end{equation*}
on the boundary. By symmetry, it holds that $(\nabla\times)^3(B\times u)_2 = 0$. The lemma is proved.

It follows from  these lemmas that
\begin{thm}\label{uniform} Let $u_0\in W\cap H^3(\Omega)$, $B_0\in \hat{W}\cap H^3(\Omega)$. Then there is a $T_0>0$ depending only on $\|(u_0,B_0)\|_{H^3}$ such that the strong solution $u=u(\nu,\mu)$, $B=B(\nu,\mu)$ of the MHD system $\ef{1.1}-\ef{1.3}$ with the initial data $(u_0,B_0)$ has the following uniform bound
\begin{equation}\label{a5.7}
\|u(\cdot,t)\|_3^2+\|B(\cdot,t)\|_3^2 + \int_0^t\nu\|u(\cdot,t)\|_4^2+\mu\|B(\cdot,t)\|_4^2\leq C\ \ \text{for} \ t\in [0,T_0],
\end{equation}
where $C$ is a constant independent of $\nu$ and $\mu$.
\end{thm}
Before a rigorous proof, we first present a formal proof for smooth solutions to \ef{1.1}-\ef{1.3} which satisfy
\begin{equation}
\begin{split}
&\partial_t(-\triangle u)-\nu\triangle(-\triangle u)-\triangle(\nabla\times u\times u+B\times\nabla\times B)=0 \ \ \text{in} \ \Omega,\\
&\partial_t(\triangle B)-\mu\triangle^2B+\triangle(B\cdot \nabla u-u\cdot\nabla B)=0 \ \ \text{in} \ \Omega,\\
&\nabla\cdot\triangle u=0,\ \ \nabla\cdot\triangle B=0 \ \ \text{in} \ \Omega,\\
&\triangle u\cdot n=0,\ \ (\nabla\times\triangle u)\times n=0 \ \ \text{on} \ \partial \Omega,\\
&\ \ \triangle B\times n=0 \ \ \text{on} \ \partial \Omega.
\end{split}
\end{equation}
Indeed, it follows from the equation of $B$ in (\ref{1.1}) with the boundary condition $B\times n=0$ and the Lemma \ref{bc1} that
\begin{equation}\label{bn}
 \triangle B\times n=0 \ \ \text{on} \ \partial \Omega.
\end{equation}
Taking the curl of the equation of $u$ in (\ref{1.1}) and using the Lemma \ref{bc2} and the boundary condition of $u$, one can get
\begin{equation}\label{un}
(\nabla\times)^3 u\times n=-(\nabla\times\triangle u)\times n=0 \ \ \text{on} \ \partial \Omega.
\end{equation}
Thus
\begin{equation}\label{5.9}
\begin{split}
\fr{d}{dt}&(\|\nabla\times(-\triangle u)\|^2+\|\nabla\times(-\triangle B)\|^2)+2(\nu\|\triangle^2u\|^2+\mu\|\triangle^2B\|^2)\\
&+2(-\triangle H_1(u,B),\triangle^2u)+2(\triangle H_2(u,B),
-\triangle^2B)=0.
\end{split}
\end{equation}
We claim that $T^*(\nu,\mu)$ is bounded below for all $\nu,\mu>0$.

Due to the boundary condition, one can integrate by part to obtain that
\begin{equation}\label{5.2a}
(-\triangle H_1(u,B),\triangle^2u)=((\nabla\times)^3H_1(u,B),\nabla\times(-\triangle u)),
\end{equation}
and
\begin{equation}\label{5.2b}
(\triangle H_2(u,B),-\triangle^2B)=((\nabla\times)^3H_2(u,B),\nabla\times(-\triangle B)),
\end{equation}
here Lemma \ref{bc3} has been used. It remains to estimate $\ef{5.2a}$ and $\ef{5.2b}$.

Noting that
\begin{equation*}
\begin{split}
&(-\triangle H_1(u,B),\triangle^2u)=((\nabla\times)^3H_1(u,B),\nabla\times(-\triangle u))\\
&=((\nabla\times)^3(B\times(\nabla\times B)),(\nabla\times)^3u) + ((\nabla\times)^3(\nabla\times u\times u),(\nabla\times)^3u)
\end{split}
\end{equation*}
and
\begin{equation*}
\begin{split}
&(\triangle H_2(u,B),-\triangle^2B)=((\nabla\times)^3H_2(u,B),\nabla\times(-\triangle B)),\\
\end{split}
\end{equation*}
one can calculate that
\begin{equation*}
 ((\nabla\times)^3(B\times(\nabla\times B)),(\nabla\times)^3u)
=-(B\cdot\nabla((\nabla\times)^3B),(\nabla\times)^3u) + R_1
\end{equation*}
where $R_1$ can be estimated by the $H^3$ norm so that
\begin{equation*}
|R_1|\leq C(\|u\|_3^3+\|B\|_3^3),
\end{equation*}
\begin{equation*}
((\nabla\times)^3(\nabla\times u\times u),(\nabla\times)^3u)
=(u\cdot\nabla(\nabla\times)^3u,(\nabla\times)^3u)+R_2
\end{equation*}
for some $|R_2|\leq C\|u\|_3^3$.
On the other hand,
\begin{equation*}
\begin{split}
& ((\nabla\times)^4(B\times u),(\nabla\times)^3B)\\
& = ((\nabla\times)^3(u\cdot\nabla B - B\cdot\nabla u),(\nabla\times)^3B) \\
& = (u\cdot\nabla (\nabla\times)^3B - B\cdot\nabla(\nabla\times)^3u,(\nabla\times)^3B) + R_3 \\
& = -(B\cdot\nabla(\nabla\times)^3u,(\nabla\times)^3B) + R_3 \\
\end{split}
\end{equation*}
for some $|R_3|\leq C(\|u\|_3^3+\|B\|_3^3)$.
Since
\begin{equation}\label{*}
\begin{split}
&(B\cdot\nabla(\nabla\times)^3u,(\nabla\times)^3B)+(B\cdot\nabla(\nabla\times)^3B,(\nabla\times)^3u) \\
& = (B,\nabla((\nabla\times)^3u\cdot(\nabla\times)^3B))  \\
& = \int_{\partial\Omega}((\nabla\times)^3u\cdot(\nabla\times)^3B)B\cdot n
- \int_\Omega((\nabla\times)^3u\cdot(\nabla\times)^3B)\nabla\cdot B=0,
\end{split}
\end{equation}
here we have used
\begin{equation*}
(\nabla\times)^3u\times n= 0,\ \  (\nabla\times)^3B\cdot n= 0, \ \ \text{on} \ \ \partial\Omega
\end{equation*}
from \ef{bn} and \ef{un}. Then we conclude that
\begin{equation*}
\begin{split}
&\frac{d}{dt}(\|(\nabla\times)^3u\|^{2}+\|(\nabla\times)^3B\|^{2})+\nu\|\Delta^2 u\|^{2}+\mu\|\Delta^2 B\|^{2} \\
&\leq C(\|u\|_3^{3}+\|B\|_3^{3}),
\end{split}
\end{equation*}
Combining it with the energy inequality (\ref{3.1}) yields that
\begin{equation*}
\begin{split}
&\frac{d}{dt}(\|u\|_3^{2}+\|B\|_3^{2})+\nu\|u\|_4^{2}+\mu\|B\|_4^{2} \\
&\leq C(\|u\|_3^{3}+\|B\|_3^{3}),
\end{split}
\end{equation*}
where $C$ is independent of $\nu$ and $\mu$, and the norm $\|\cdot\|_s$ is the equivalent one in (\ref{2.3}).
Comparing with the ordinary differential equation
\begin{equation*}
y'(t)=Cy(t)^{\fr{3}{2}},
\end{equation*}
\begin{equation*}
y(0)=\|u(0)\|_3^2+\|B(0)\|_3^2,
\end{equation*}
and let $T_0$ be the blow up time, one obtains that
\begin{equation*}
T^*(\nu,\mu)\geq T_0\quad \text{for all} \ \nu, \mu>0,
\end{equation*}
and (\ref{a5.7}) is valid.\\
\textbf{Proof of Theorem \ref{uniform}:}
To make the proof rigorous, one can use the Galerkin approximations. Consider the system satisfied by  $-\triangle u^m(x,t)$ and $-\triangle B^m(x,t)$. Let $u_0\in W$, $B_0\in \hat{W}$. It follows from $\ef{3.5}$ and $\ef{3.6}$ that
\begin{equation}\label{5.13}
(-\triangle u^m)'-\nu\triangle(-\triangle u^m)+\sum g^1_j\lambda_je_j=0,
\end{equation}
\begin{equation}\label{5.14}
(-\triangle B^m)'-\mu\triangle(-\triangle B^m)+\sum g^2_j\lambda_j\hat {e}_j=0,
\end{equation}
\begin{equation}\label{5.15}
(-\triangle u^m)(0)=P_m(-\triangle u_0),\quad (-\triangle B^m)(0)=\hat{P}_m(-\triangle B_0).
\end{equation}
Since $\nabla\times e_i\times n=0$ and $n\times\nabla\times H_1(u^m,B^m)=0$ on the boundary, so integration by parts yields
\begin{equation*}
\begin{split}
&(-\triangle P_mH_1(u^m,B^m),e_i)=(\sum g^1_j\lambda_j e_j,e_i)\\
&=(H_1(u^m,B^m),-\triangle e_i)\\
&=\int_{\partial\Omega}H_1(u^m,B^m)\cdot(\nabla\times e_i\times n)+(\nabla\times H_1(u^m,B^m),\nabla\times e_i)\\
&=\int_{\partial\Omega}(n\times\nabla\times H_1(u^m,B^m))\cdot e_i+(-\Delta H_1(u^m,B^m), e_i)\\
&=(-\triangle H_1(u^m,B^m),e_i).
\end{split}
\end{equation*}

Thus the  following  commutation holds
\[   \Delta P_mH_1(u^m,B^m) =  P_m\Delta H_1(u^m,B^m)  \]
where
\[   P_mH_1(u^m,B^m)=\sum g^1_je_j   \]
with $g^1_j = (H_1(u^m,B^m),e_j)$, and
\[   P_m\Delta H_1(u^m,B^m)=\sum g^{1,\delta}_je_j   \]
with $g^{1,\delta}_j = (\Delta H_1(u^m,B^m),e_j)$.

Similarly, integration by parts shows that
\begin{equation*}
\begin{split}
&(-\Delta \hat{P}_mH_2(u^m,B^m),\hat{e}_i)=(\sum g^2_j\lambda_j \hat{e}_j,\hat{e}_i)\\
&=(H_2(u^m,B^m),-\triangle \hat{e}_i)\\
&=\int_{\partial\Omega}(n\times H_2(u^m,B^m))\cdot(\nabla\times\hat{e}_i)+(\nabla\times H_2(u^m,B^m),\nabla\times \hat{e}_i)\\
&=\int_{\partial\Omega}\nabla\times H_2(u^m,B^m)\cdot (\hat {e}_i\times n)+(-\Delta H_2(u^m,B^m), \hat{e}_i)\\
&=(-\triangle H_2(u^m,B^m),\hat{e}_i)
\end{split}
\end{equation*}
due to the fact that $\hat{e}_i\times n=0$ and $n\times H_2(u^m,B^m)=0$ on the boundary.
Thus the following commutation holds
\[   \Delta \hat{P}_mH_2(u^m,B^m) =  \hat{P}_m\Delta H_2(u^m,B^m)  \]
where
\[   \hat{P}_mH_2(u^m,B^m)=\sum g^2_j\hat {e}_j   \]
with $g^2_j = (H_2(u^m,B^m),\hat{e}_j)$, and
\[   \hat{P}_m\Delta H_2(u^m,B^m)=\sum g^{2,\delta}_j\hat{e}_j   \]
with $g^{2,\delta}_j = (\Delta H_2(u^m,B^m),\hat {e}_j)$.
It follows from (\ref{5.13}), (\ref{5.14}) and the above commutations that
\begin{equation}\label{5.16}
\begin{split}
\fr{d}{dt}&(\|\nabla\times(-\triangle u^m)\|^2+\|\nabla\times(-\triangle B^m)\|^2)+2(\nu\|\triangle^2u^m\|^2+\mu\|\triangle^2B^m\|^2)\\
&+2(\triangle H_1(u^m,B^m),-\triangle^2u^m)+2(\triangle H_2(u^m,B^m),-\triangle^2B^m)=0.
\end{split}
\end{equation}
Note that $(\nabla\times\Delta u^m)\times n = 0$ on the boundary. It follows that
\[    (\triangle H_1(u^m,B^m),-\triangle^2u^m) = ((\nabla\times)^3 H_1(u^m,B^m),(\nabla\times)^3u^m).   \]
By Lemma \ref{bc3}, it follows that $(\triangle H_2(u^m,B^m))\times n = 0$ on the boundary, and then
\[   (\triangle H_2(u^m,B^m),-\triangle^2B^m) = ((\nabla\times)^3H_2(u^m,B^m),(\nabla\times)^3B^m).     \]
Hence, the estimates in the formal analysis above can also be applied to the Galerkin approximations, and the corresponding bounds can also be obtained, which allow one to pass the limit to derive the desired a priori estimates \ef{a5.7}. The theorem is proved.

The above uniform estimates allow us to obtain the zero dissipation limit.
\begin{thm}\label{h3}  Assume that $u_0\in W\cap H^3(\Omega)$, $B_0\in \hat{W}\cap H^3(\Omega)$. Let $(u,B)=(u(\nu,\mu),B(\nu,\mu))$ be the corresponding strong solution to the MHD system $\ef{1.1}-\ef{1.3}$ on $[0,T_0]$ in Theorem \ref{uniform}. Then as $\nu,\mu\longrightarrow 0$, $(u,B)$ converges to the unique solution $(u^0,B^0)$ of the ideal MHD system with the same initial data in the sense
\begin{equation}\label{5.17}
u(\nu,\mu), B(\nu,\mu)\longrightarrow u^0, H^0\ \text{in}\ L^q(0,T_0;H^2(\Omega)).
\end{equation}
\begin{equation}\label{5.18}
u(\nu,\mu), B(\nu,\mu)\longrightarrow u^0, H^0\ \text{in}\ C(0,T_0;H^2(\Omega))
\end{equation}
for all $1\leq q<\infty$.
\end{thm}
\pf It follows from theorem \ref{uniform} that
\begin{equation*}
u(\nu,\mu), B(\nu,\mu) \ \text{is uniformly bounded in }\ C([0,T_0];H^3(\Omega)),
\end{equation*}
\begin{equation*}
u'(\nu,\mu), B'(\nu,\mu) \ \text{is uniformly bounded in }\ L^2(0,T_0;H^2(\Omega)),
\end{equation*}
for all $\nu, \mu>0$. From the Aubin-Lions lemma, there is a subsequence $\nu_n$, $\mu_n$ and $u^0,B^0$ such that
\begin{equation*}\begin{split}
&(u(\nu_n,\mu_n),B(\nu_n,\mu_n))\longrightarrow(u^0,B^0) \ \text{in}\ L^\infty(0,T_0;H^3(\Omega))\ \ \text{weakly},\\
&(u(\nu_n,\mu_n),B(\nu_n,\mu_n))\longrightarrow(u^0,B^0) \ \text{in}\ L^p(0,T_0;H^2(\Omega)),\\
&(u(\nu_n,\mu_n),B(\nu_n,\mu_n))\longrightarrow (u^0,B^0) \ \text{in}\ C([0,T_0);H^2(\Omega))
\end{split}
\end{equation*}
for any $1\leq p<\infty$ as $\nu_n$,  $\mu_n\longrightarrow0$. Passing to the limit shows that $(u^0,B^0)$ solves the following limit equations
\begin{equation}\label{imhd}
\begin{split}
&\partial_t{ u^0}+(\nabla\times u^0)\times u^0+B^0\times \nabla\times B^0+\nabla p=0\ \ \text{in} \ \  \Omega,\\
&\nabla\cdot u^0=0\ \  \text{in} \ \ \Omega,\\
&\partial_t{B^0}=\nabla\times (u^0\times B^0)\ \text{in} \ \  \Omega,\\
&\nabla\cdot B^0=0\ \ \text{in} \ \Omega
\end{split}
\end{equation}
with the boundary conditions
\begin{equation}\label{ibc}\begin{split}
&u^0\cdot n=0,\quad (\nabla\times u^0)\times n=0 \ \ \text{on} \ \ \partial \Omega,\\
&B^0\times n=0\ \  \text{on} \ \ \partial \Omega,\\
&\triangle B^0\times n=0 \ \  \text{on} \ \ \partial \Omega,
\end{split}
\end{equation}
and $p$ satisfying
\begin{equation}\label{ip}
\begin{split}
&\triangle p=-\big(\nabla\cdot((\nabla\times u^0)\times u^0)-\nabla\cdot((\nabla\times B^0)\times B^0)\big),\\
&\nabla p\cdot n=0 \ \  \text{on} \ \ \partial \Omega.
\end{split}
\end{equation}
Let $(u,B)$ and $(u^0,B^0)$ be two strong solutions to \ef{imhd}-\ef{ip}. Set $\bar{u}=u-u^0$, and $\bar{B}=B-B^0$. Then
\begin{equation}\label{4.8a}
\partial_t\bar{u}+H_1(u,B)-H_1(u^0,B^0)=0,
\end{equation}
\begin{equation}\label{4.9a}
\partial_t\bar{B}+H_2(u,B)-H_2(u^0,B^0)=0.
\end{equation}
Taking the inner products with $\bar{u}$ in $\ef{4.8a}$, and  $\bar{B}$ in $\ef{4.9a}$ and integrating by parts lead to
\begin{equation*}
\begin{split}
\fr{d}{dt}(\|\bar{u}\|^2+\|\bar{B}\|^2)
+(H_1(u,B)-H_1(u^0,B^0),\bar{u})+(H_2(u,B)-H_2(u^0,B^0),\bar{B})=0.\\
\end{split}
\end{equation*}
\ef{h1h2} implies that
\begin{equation}\label{ih1h2}
\begin{split}
&(H_1(u,B)-H_1(u^0,B^0),\bar{u})+(H_2(u,B)-H_2(u^0,B^0),\bar{B})\\
&=(\nabla\times\bar{u}\times u,\bar{u})+(\nabla\times(\bar{u}\times\bar{B}),B^0)+(\bar{B}\times u^0,\nabla\times \bar{B})).\\
\end{split}
\end{equation}
Since
\begin{equation}\label{iu}
(\nabla\times\bar{u}\times u,\bar{u})=((\bar{u}\cdot\nabla) u,\bar{u})-((u\cdot \nabla) \bar{u},\bar{u})=((\bar{u}\cdot\nabla) u,\bar{u}),
\end{equation}
\begin{equation}\label{iub}
(\nabla\times(\bar{u}\times\bar{B}),B^0)=(\bar{u}\times\bar{B},\nabla\times B^0),
\end{equation}
and
\begin{equation}\label{iB}
(\bar{B}\times u^0,\nabla\times \bar{B})=((u^0\cdot\nabla)\bar{B},\bar{B})-((\bar{B}\cdot\nabla)u^0,\bar{B})=((u^0\cdot\nabla)\bar{B},\bar{B}),
\end{equation}
here the boundary conditions $B\times n=0$, $B^0\times n=0$, $u\cdot n=0$ and $u^0\cdot n=0$ have been used,
one can get
\begin{equation*}
\begin{split}
\fr{d}{dt}(\|\bar{u}\|^2+\|\bar{B}\|^2)
\leq C(\|u^0\|_3+\|B^0\|_3+\|u\|_3)(\|\bar{u}\|^2+\|\bar{B}\|^2).
\end{split}
\end{equation*}
Note that $u^0$, $B^0$, $u$ and $B$ are all in $L^\infty(0,T_0;H^3)$ and $\bar{u}(0)=\bar{B}(0)=0$. One obtains the uniqueness by Gronwall's inequality.

Finally, we prove the following convergence rate.
\begin{thm}\label{h4} Under the same assumptions in Theorem \ref{h3}, it holds that
\begin{equation*}
\|u(\nu,\mu)-u^0\|^2_2+\|B(\nu,\mu)-B^0\|^2_2\leq C(T_0)(\nu +\mu)
\end{equation*}
on the interval $[0,T_0]$.
\end{thm}

\pf Set $\bar{u}=u(\nu,\mu)-u^0$ and $\bar{B}=B(\nu,\mu)-B^0$. One can get that $-\triangle\bar{u}$ and $-\triangle\bar{B}$ solve
\begin{equation}\label{5.6}
\partial_t(-\triangle\bar{u})-\triangle(H_1(u,B)-H_1(u^0,B^0))=-\nu\triangle^2u\ \ \text{in}\ \ \Omega,
\end{equation}
\begin{equation}\label{5.7}
\partial_t(-\triangle\bar{B})-\triangle(H_2(u,B)-H_2(u^0,B^0))=-\mu\triangle^2B\ \ \text{in}\ \ \Omega,
\end{equation}
\begin{equation}\label{5.8}
\nabla\cdot{\bar{u}}=0,\quad \nabla\cdot{\bar{B}}=0 \ \ \text{in}\ \ \Omega,
\end{equation}
\begin{equation}\label{5.9}
\bar{u}\cdot n=0,\quad \bar{B}\times n=0, \ \ \text{on}\ \ \partial \Omega,
\end{equation}
with $\nabla\times u\times n=0$, $\nabla\times u^0\times n=0$, $\triangle B\times  n=0$, $(\nabla\times)^3 u\times  n=0$, and $ \triangle B^0\times  n=0$ on the boundary. Taking inner product of $\ef{5.6}$ with $-\triangle\bar{u}$ and $\ef{5.7}$ with
 $-\triangle\bar{B}$ and integrating by parts, one gets that
 \begin{equation*}
 \begin{split}
 &\fr{d}{dt}(\|\triangle\bar{u}\|^2+\|\triangle\bar{B}\|^2)-2(\triangle(H_1(u,B)-H_1(u^0,B^0)),-\triangle\bar{u})\\
 &-2(\triangle(H_2(u,B)-H_2(u^0,B^0)),-\triangle\bar{B})\\
 &=\nu((\nabla\times)^3u,(\nabla\times)^3\bar{u})+\mu((\nabla\times)^3B,(\nabla\times)^3\bar{B}).
 \end{split}
 \end{equation*}
A simple computation yields
\begin{equation*}
\begin{split}
-\triangle(H_1(u,B)-H_1(u^0,B^0))&=(u\cdot\nabla)(-\triangle\bar{u})-(B\cdot\nabla)(-\triangle\bar{B})\\
&+(\bar{u}\cdot\nabla)(-\triangle u^0)-(\bar{B}\cdot\nabla)(-\triangle B^0)\\
&+\sum_{i,j=1,2,i+j=3}F_{i,j}(D^iu^0,D^j\bar{u})-\sum_{i,j=1,2,i+j=3}F_{i,j}(D^iB^0,D^j\bar{B})\\
&+\sum_{i,j=1,2,i+j=3}F_{i,j}(D^iu,D^j\bar{u})-\sum_{i,j=1,2,i+j=3}F_{i,j}(D^iB,D^j\bar{B}),
\end{split}
\end{equation*}
\begin{equation*}
\begin{split}
-\triangle(H_2(u,B)-H_2(u^0,B^0))&=(u\cdot\nabla)(-\triangle\bar{B})-(B\cdot\nabla)(-\triangle\bar{u})\\
&+(\bar{u}\cdot\nabla)(-\triangle B^0)-(\bar{B}\cdot\nabla)(-\triangle u^0)\\
&+\sum_{i,j=1,2,i+j=3}F_{i,j}(D^iu^0,D^j\bar{B})
-\sum_{i,j=1,2,i+j=3}F_{i,j}(D^iB^0,D^j\bar{u})\\
&+\sum_{i,j=1,2,i+j=3}F_{i,j}(D^iu,D^j\bar{B})-\sum_{i,j=1,2,i+j=3}F_{i,j}(D^iB,D^j\bar{u}),
\end{split}
\end{equation*}
where $F_{i,j}(D^iu,D^jv)'s$ are bilinear forms and ${D^i}'s$ are the i-th order differential operators.
It follows from $\nabla\cdot u =0$ and $u\cdot n|_{\partial\Omega}=0$ that
\begin{equation*}
\big((u\cdot\nabla)(-\triangle\bar{u}),-\triangle \bar{u}\big)=0, \ \ \big((u\cdot\nabla)(-\triangle\bar{B}),-\triangle \bar{B}\big)=0.
\end{equation*}
On the other hand,
\begin{equation}\label{**}
\begin{split}
&\big((B\cdot\nabla)(-\triangle \bar{B}),-\triangle \bar{u}\big)+\big((B\cdot\nabla)(-\triangle\bar{u}),-\triangle\bar{B}\big) \\
& = (B,\nabla(\triangle \bar{B}\cdot\triangle\bar{u}))  \\
& = \int_{\partial\Omega}(\triangle \bar{B}\cdot\triangle\bar{u})B\cdot n
- \int_\Omega(\triangle \bar{B}\cdot\triangle\bar{u})\nabla\cdot B=0,
\end{split}
\end{equation}
where one has used
\begin{equation*}
\triangle \bar{B}\times n =0,\ \ \triangle \bar{u}\cdot n=0,\ \ \text{on} \ \ \partial\Omega.
\end{equation*}
 Therefore,
\begin{equation*}
\begin{split}
&|(\triangle(H_1(u,B)-H_1(u^0,B^0)),-\triangle\bar{u})+(\triangle(H_2(u,B)-H_2(u^0,B^0)),-\triangle\bar{B})|\\
&\leq C(\|u^0\|_3+\|B^0\|_3+\|u\|_3+\|B\|_3)(\|\triangle\bar{u}\|^2+\|\triangle\bar{B}\|^2).
\end{split}
\end{equation*}
Also, one has that
\begin{equation*}
|((\nabla\times)^3u,(\nabla\times)^3\bar{u})|\leq C\|(\nabla\times)^3u\|(\|(\nabla\times)^3u\|+\|(\nabla\times)^3u^0\|)
\end{equation*}
and
\begin{equation*}
|((\nabla\times)^3B,(\nabla\times)^3\bar{B})|\leq C\|(\nabla\times)^3B\|(\|(\nabla\times)^3B\|+\|(\nabla\times)^3B^0\|).
\end{equation*}
These estimates are uniform respect to $\nu$, $\mu$ and thus
\begin{equation*}
\fr{d}{dt}(\|\triangle\bar{u}\|^2+\triangle\bar{B}\|^2)\leq C(T_0)(\|\triangle\bar{u}\|^2+\|\triangle\bar{B}\|^2+\nu+\mu).
\end{equation*}
Due to $\bar{u}(0)=0$,  $\bar{B}(0)=0$ and Gronwall's inequality, we deduce that
\begin{equation}\label{5.10}
\|\triangle\bar{u}\|^2+\|\triangle\bar{B}\|^2\leq C(T_0)(\nu +\mu).
\end{equation}
The theorem is proved.

{\bf Acknowledge.}
 Yuelong Xiao is partially supported by NSFC Nos. {\rm 11871412, 11771300}. And  Qin Duan is partially supported by NSFC Nos. {\rm 11771300}. The research of Zhouping Xin was supported in part by Zheng Ge Ru Foundation, HongKong RGC Earmarked Research Grants: CUHK14302819, CUHK14300917,  CUHK14302917, and Basic and Applied Basic Research Foundation of Guangdong Province 2020B1515310002.


\begin{thebibliography}{29}
	\bibitem{BF} H. Beir\~ao da Veiga and F. Crispo: Sharp inviscid limit results under
	Navier type boundary conditions. An Lp theory, J. math. fluid mech.
	12 (2010), 397-411.
	\bibitem{BF1} H. Beir\~ao da Veiga and  F. Crispo:  Concerning the $W^{k,p}$-inviscid limit for 3D
	flows under a slip boundary condition, J. math. fluid mech.
	13 (2011), 117-135.
	\bibitem{BF2} H. Beir\~ao da Veiga and  F. Crispo: The 3D inviscid limit result under slip boundary conditions. A negative
	answer J. Math. Fluid Mech. 14(2012) 55-59.
	\bibitem{BF3}  H. Beir\~ao da Veiga and  F. Crispo: A missed persistence property for the Euler equations
	and its effect on inviscid limits, Nonlinearity 25 (2012), 1661-1669.
	\bibitem{BN08} Bellout, H.and Neustupa, J.: A Navier-Stokes approximation of
	the 3D Euler equation with the zero flux on the boundary.
	\bibitem{BN10} Bellout, H. and Neustupa, J. and Penel, P.: On a $\nu$ continous family of strong
	solution to the Euler or Navier-Stokes equations with the Navier type boundary condition, Disc. Cont. Dyn. Sys. Vol. 27:4 (2010), 1353-1373.
	\bibitem{Ber}Berselli, Luigi C.: Some results on the Navier-Stokes equations
	with Navier boundary conditions. Riv. Math. Univ. Parma (N.S.) 1 (2010), no. 1, 1-75.
	\bibitem{Ber1} Berselli, Luigi C. and Spirito, S.: On the Vanishing Viscosity Limit of 3D Navier-Stokes
	Equations under Slip Boundary Conditions in General Domains, Commun. Math. Phys. (DOI) 10.1007.
	\bibitem{CMR} Clopeau,T. and Mikeli\'{c},A. and Robert,R.: On the
	valishing viscosity limit for the 2D incompressible Navier-Stokes
	equations with the friction type boundary conditions,
	Nonlinearity 11(1998), 1625-1636.
	\bibitem{Co1} Constantin, P.: Note on loss of regularity for solutions of the 3-D incompressible Euler and related
	equations. Commun. Math. Phys. 104, 311-326 (1986).
	\bibitem{Co}Constantin,P. and
	Foias,C.:  Navier Stokes equation, Univ. of Chicago press
	IL(1988).
	\bibitem{DL} Duvaut, G. and Lions, J.-L.: In\'{e}quations en thermo\'{e}lasticit\'{e} et magn\'{e}tohydrodynamique, Arch. Rational Mech. Anal.,
46, 1972, 241-279.
	\bibitem{Gun} Gunzburger, Max D. and Meir, Amnon J. and Peterson, Janet S.: On the existence, uniqueness, and finite element approximation of solutions of the equations of stationary, incompressible magnetohydrodynamics. Math. Comp. 56 (1991), no. 194, 523-563.
	\bibitem{HHM} Han, Woo Jin and Hwang, Hyung Ju and Moon, Byung Soo: On the well-posedness of the Hall-magnetohydrodynamics with the ion-slip effect, (English summary) J. Math. Fluid Mech., 21 (2019), no. 4, Paper No. 47, 28 pp.
	\bibitem{HY} Hughes, William and Young, F.: The electromagnetodynamics of fluids, SERBIULA (sistema Librum 2.0), 05, 2020.
		\bibitem{IfP} Iftimie, D. and Planas, G.: Inviscid limits for the
	Navier-Stokes equations with Navier friction boundary conditions. Nonlinearity 19 (2006), no. 4, 899-918.
	\bibitem{IfS}  Iftimie, D. and Sueur, F.: Viscous boundary layers for the Navier-Stokes
	equations with the Navier slip conditions. Arch. Ration. Mech. Anal. 199 (2011), no. 1, 145-175.
	\bibitem{Ja} Jackson, John David:
	Classical electrodynamics, Second, John Wiley \& Sons, Inc., New York-London-Sydney, 1975.
	\bibitem{JeanF06} Gerbeau, Jean-Fr¨¦d¨¦ric and Le Bris, Claude and Leli¨¨vre, Tony: Mathematical methods for the magnetohydrodynamics of liquid metals. Numerical Mathematics and Scientific Computation. Oxford University Press, Oxford, 2006.
	\bibitem{Ka} Kato, T.: Nonstationary flows of viscous and ideal fluids in R3.
	J. Functional Analysis 9 (1972), 296-305.
	\bibitem{Ka1}Kato, T.: Quasi-linear equations of evolution, with applications
	to partial differential equations. Spectral theory and differential equations
	(Proc. Sympos., Dundee, 1974; dedicated to Konrad J$\ddot{o}$rgens), pp. 25-70.
	Lecture Notes in Math., Vol. 448, Springer, Berlin, 1975.
	\bibitem{Ka2} Kato,T.:
	Remarks on zero viscosity limit for non-stationary
	Navier-Stokes flows with boundary. In: Seminar on PDE (S.S.Chen,
	eds), Springer, New York, (1984), 85-98.
	\bibitem{KaL} Kato, T. and Lai, C. Y.: Nonlinear evolution equations and the Euler flow.
	J. Funct. Anal. 56 (1984), no. 1, 15-28.
	\bibitem{Ke1} Kelliher, J.: Navier-Stokes equations with Navier boundary conditions
	for a bounded domain in the plane. SIAM J. Math. Anal. 38 (2006), no. 1, 210-232 (electronic).
	\bibitem{Ke2} GM Gie and JP Kelliher: Boundary layer analysis of the Navier¨CStokes equations with generalized
	Navier boundary conditions, JDE 253:6, 1862-1892.
	\bibitem{LXY} Liu, Cheng-Jie and Xie, Feng and Yang, Tong: MHD boundary layers theory in Sobolev spaces without monotonicity I: Well-posedness theory. Comm. Pure Appl. Math. 72 (2019), no. 1, 63-121.
	\bibitem{Mae} Maekawa, Yasunori: On the inviscid limit problem of the vorticity equations for viscous incompressible flows in the half-plane. Comm. Pure Appl. Math. 67 (2014), no. 7, 1045-1128.
	\bibitem{Mas} Masmoudi, N.: Remarks about the inviscid limit of the Navier-Stokes system.
	Comm. Math. Phys. 270 (2007), no. 3, 777-788.
	\bibitem{Mas1} Masmoudi, N. and Rousset, F.: Uniform Regularity for the Navier-Stokes equation
	with Navier boundary condition, Arch. Rational Mech. Anal. 203(2012), 529-575.
	\bibitem{PLL96} Lions, P.-L.: Mathematical topics in fluid mechanics. Vol. 1. Incompressible models.
	Oxford Lecture Series in Mathematics and Its Applications, 3. Oxford Science
	Publications. The Clarendon Press, Oxford University Press, New York, 1996.
	\bibitem{she} Shercliff, J. A.:
	A textbook of magnetohydrodyamics, Pergamon Press, Oxford-New York-Paris, 1965.
	\bibitem{Sm1} Sammartino, M. and Caflisch, R. E.: Zero viscosity limit for analytic solutions of the Navier-Stokes
	equation on a half-space. I. Existence for Euler and Prandtl equations. Comm. Math. Phys.
	192 (1998), no. 2, 433-461.
	\bibitem{Sm2} Sammartino, M.and Caflisch, R. E.: Zero viscosity limit for analytic solutions of the Navier-Stokes
	equation on a half-space. II. Construction of the Navier-Stokes solution. Comm. Math. Phys.
	192 (1998), no. 2, 463-491.
	\bibitem{ST}Sermange, Michel and Temam, Roger:
Some mathematical questions related to the {MHD} equations, Comm. Pure Appl. Math.,
36, 1983, 5, 635-664.
	\bibitem{StoG} Str$\ddot{o}$hmer, Gerhard: About an initial-boundary value problem from magnetohydrodynamics. Math. Z. 209 (1992), no. 3, 345-362.
	\bibitem{Wa} Wang, X.: A Kato type theorem on zero viscosity limit of Navier-Stokes flows.
	Dedicated to Professors Ciprian Foias and Roger Temam (Bloomington, IN, 2000).
	Indiana Univ. Math. J. 50 (2001), Special Issue, 223-241.
	\bibitem{WWX} Wang, X. and  Wang, Y.and Xin, Z.: Boundary layers in incompressible Navier-Stokes
	equations with Navier boundary conditions for the vanishing viscosity limit.
	Commun. Math. Sci. 8 (2010), no. 4, 965-998.
	\bibitem{WXZ} Wang, L. and Xin, Z. and Zang, A.: Vanishing Viscous Limits for 3D Navier-Stokes Equations
	with A Navier-Slip Boundary Condition. Arxiv preprint arXiv:1201.1986, 2012 - arxiv.org.
	\bibitem{XX} Xiao, Y.L. and Xin, Z.P.: On the vanishing viscosity limit for the 3D
	Navier-Stokes equations with a slip boundary condition. Comm. Pure Appl. Math. Vol.
	LX (2007) 1027-1055.
	\bibitem{XX1}Xiao, Y.L. and Xin, Z.P.: Remarks on the vanishing viscosity limit for 3D
	Navier-Stokes equations with a slip boundary condition, Chinese
	Ann. Math., 32B(3)(2011), 321-332.
	\bibitem{XX2}Xiao, Y. and Xin, Z.: On 3D Lagrangian Navier-Stokes $\alpha$
	model with a Class of Vorticity Slip Boundary conditions. J. Math. Fluid Mech. 15 (2013), 215-247.
	\bibitem{XX3}Xiao, Y.L. and Xin, Z.P.: A New Boundary Condition for the 3D Navier-Stokes Equation and
	the Vanishing Viscosity Limit, J. Math. Phys. 53, 115617 (2012).
	\bibitem{XX4}Xiao, Y.L. and Xin, Z.P.: On the inviscid limit of the 3D Navier-Stokes equations with generalized Navier-slip boundary conditions. Commun. Math. Stat. 1 (2013), no. 3, 259-279.
	\bibitem{XXW}Xiao, Y.L. and Xin, Z.P. and Wu, J.H.: Vanishing viscosity limit for the 3D magnetohydrodynamic system with a slip boundary condition. J. Funct. Anal. 257 (2009), no. 11, 3375-3394.
	
\end{thebibliography}
\end{document}